\documentclass[preprint,12pt]{elsarticle}
\journal{Computers and Structures}
\usepackage{graphics,epsfig}
\usepackage{graphicx}
\usepackage{float}
\usepackage{amssymb,amstext,amsmath}
\usepackage{mathtools,physics,verbatim}
\usepackage{booktabs}
\usepackage{natbib}
\usepackage{xspace}
\usepackage{lscape}
\usepackage{multirow}
\usepackage[htt]{hyphenat}
\usepackage{xcolor}

\newcommand{\gb}[1]{\textcolor{black}{#1}} 

\begin{document}
\begin{frontmatter}
\title{Asymptotically accurate and locking-free finite element implementation of first order shear deformation theory for plates}
\author{K.~C. Le$^{a,b}$\footnote{Corresponding author. Phone: +84 93 4152458, email: lekhanhchau@tdtu.edu.vn}, H.-G. Bui$^c$}
\address{$^a$Institute for Advanced Study in Technology, Ton Duc Thang University, Ho Chi Minh City, Vietnam
\\
$^b$Faculty of Civil Engineering, Ton Duc Thang University, Ho Chi Minh City, Vietnam
\\
$^c$Institute of Material Systems Modeling, Helmholtz-Zentrum Hereon, Geesthacht, Germany
\\
\bigskip
Dedicated to the 80th birthday of V.L. Berdichevsky}

\begin{abstract} 
A formulation of the asymptotically exact first-order shear deformation theory for linear-elastic homogeneous plates in the rescaled coordinates and rotation angles is considered. This allows the development of its asymptotically accurate and shear-locking-free finite element implementation. As applications, numerical simulations are performed for circular and rectangular plates, showing complete agreement between the analytical solution and the numerical solutions based on two-dimensional theory and three-dimensional elasticity theory.
\end{abstract}


\end{frontmatter}


\section{Introduction}
The first-order shear deformation theory (FSDT) for plates, originally proposed by Reissner \cite{reissner1945the},\footnote{See also Mindlin \cite{mindlin1951influence}, whose motivation for including transverse shear and rotatory inertia was somewhat different: his main goal was to capture both low- and high-frequency thickness vibrations of the plates.} has since attracted the attention of both theorists and practitioners alike, mainly because of its applicability to moderately thick plates and the development of numerical methods \cite{bath1986a}-\cite{belardi2021on}, but also because of the logic behind its derivation, which can be applied to other problems \cite{wang2000shear}-\cite{challamel2019a}. The first asymptotically exact version of FSDT for plates has been derived by Berdichevsky \cite{berdichevsky1979variational} using the variational-asymptotic method he himself developed. The extension of his result to laminated plates using the same method has been considered by Sutyrin \cite{sutyrin1997derivation} and Yu \cite{yu2005mathematical}. However, since in the general case of laminated plates, the dimension reduction does not lead to a FSDT, these authors have tried to optimize the parameters so that a derived theory is as close as possible to asymptotic correctness while being a FSDT. Le \cite{le2023asymptotically} has recently shown that the construction of asymptotically exact FSDT for functionally graded (FG) plates is possible when the mass density and elastic moduli vary across thickness such that their distributions are even functions of the transverse coordinate. Similar to Berdichevsky's FSDT for homogeneous plates, his theory for FG-plates is asymptotically exact up to the order of $h^2/l^2$, where $h$ is the plate thickness and $l$ is the characteristic scale of change of the deformation state in the longitudinal directions. Also worth mentioning are some recent applications of the variational-asymptotic method to dimension reduction, homogenization, nonlinear vibrations, and wave propagation in \cite{le2013}-\cite{le2020a}.

In numerical simulations of bending deformation of plates under a transverse load within FSDT based on the finite element method, the so-called shear-locking (SL) effect often occurs, especially when using low-order finite plate elements (see the above cited papers \cite{bath1986a}-\cite{belardi2021on}). Physically, this effect is due to the fact that the shear stiffness, in terms of the small plate thickness $h$, is two orders of magnitude larger than the bending stiffness, while on the other hand the rotation angles caused by pure shear are much smaller than the bending measures (or changes in curvature) of the plate, since the latter are of the order of the characteristic strain divided by the thickness. This is also evident from the variational-asymptotic analysis of the energy functional \cite{le1999vibrations,berdichevsky2009variational}. Therefore, the numerical instability with respect to the shear energy and also to the constitutive equations must occur in the limit $h\to 0$ due to the multiplication of the extremely small and large numbers when the standard low-order finite element calculation is used. There are several sophisticated methods that make it possible to alleviate this shear-locking effect but at the expense of computational efficiency. The reduced and selective integration method, using two different integration rules for bending and shear energies \cite{zienkiewicz1971reduced}-\cite{pugh1978a}, has long been the preferred technique for the numerical treatment of SL. The mathematical justification of this method, based on the equivalence between the reduced integration approach and certain mixed models, was later given by Malkus and Hughes \cite{malkus1978mixed}. Unfortunately, the reduced integration technique often leads to instability due to rank deficiency and to zero-energy modes \cite{jackson1981singular}-\cite{hayes1981ill}. Therefore, several alternative formulations and numerical techniques have been developed to mitigate SL and increase the accuracy and stability of the solution. These include the modified shear strain method \cite{hughes1981finite,crisfield1984a}, the hybrid and mixed method \cite{lee1978finite}-\cite{dolbow1999}, the extended assumed strain method \cite{simo1990a,cardoso2008enhanced}, the assumed natural strain method \cite{hughes1981finite,tessler1985a}, and the shear gap method \cite{Bletzinger2000a} (see also the recent paper \cite{nguyen2010an} discussing relevant publications). 

As far as the authors of this paper are aware, none of the existing studies on this topic has found a simple formulation of the FSDT that is inherently free of the shear-locking effect, regardless of the discretization scheme and integration technique used. The goal of this paper is therefore twofold. First, we give the formulation of the FSDT for plates in the rescaled coordinates and rotation angles. This formulation occurs naturally when the coordinates in the mid-plane are scaled by the plate thickness $h$ (thus becoming dimensionless), while the rotation angles are multiplied by $h$, resulting in equal and finite orders of the bending and shear stiffnesses as well as the scaled rotation angles and bending measures.\footnote{Compare with \cite{oesterle2018intrinsically}, where the intrinsically shear-locking-free formulation was obtained by reparametrizing the kinematic equations for bending measures and shear.}  Since this formulation is independent of the plate thickness and inherently shear-locking-free, no high-order interpolation scheme and/or sophisticated integration technique is required for the discretization and FE-implementation, so the computational efficiency can be significantly improved. However, our second goal is to ensure that the FE-implementation is asymptotically accurate. According to \cite{berdichevsky1979variational,le2023asymptotically}, the asymptotic accuracy requires that both the transverse displacement and the rotation angles belong to the $C^1$-function space. We will show by numerical simulations of the circular and rectangular plates and comparison with the analytical solutions of the FSDT and the numerical solution of the three-dimensional elasticity theory that the asymptotic accuracy is indeed achieved, provided that the isogeometric elements guaranteeing the $C^1$-continuity for the primary variables are used.

The paper is organized as follows. After this brief introduction, Section 2 gives the inherently shear-locking-free rescaled variational formulation of the FSDT for plates. Section 3 is devoted to its FE-implementation. In Section 4, we consider two numerical examples where the developed FE-code is applied: (i) circular plates subjected to uniform loading, and (ii) rectangular plates subjected to uniform loading. Finally, Section 5 concludes the paper.
 
\section{Rescaled variational formulation of FSDT for plates}
\begin{figure}[htb]
	\centering
	\includegraphics[width=8cm]{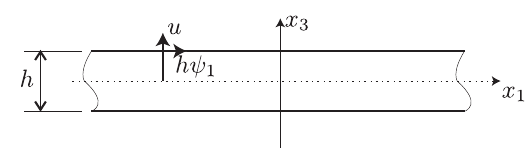}
	\caption{Cross section of a plate}
	\label{fig:1}
\end{figure}

Let $\Omega $ be a two-dimensional domain in the $(x_1,x_2)$-plane bounded by a piecewise smooth closed curve $\partial \Omega $. We consider a linear elastic homogeneous plate which in the undeformed stress-free state occupies the 3-D region $\mathcal{V}=\Omega \times (-h/2,h/2)$. Its cross section in the plane $(x_1,x_3)$ is shown in Fig.~\ref{fig:1}. We call $\Omega $ the plate mid-plane and $h$ its thickness. 

Variational principle of the asymptotically exact two-dimensional first-order shear deformation theory (FSDT) for linearly elastic homogeneous plate \cite{berdichevsky1979variational,berdichevsky2009variational} states that the true deflection and rotation angles of the plate minimize the 2-D average energy functional 
\begin{equation}\label{energyfunctional}
J[u,\varphi_\alpha ]=\int_{\Omega }\Bigl\{ \frac{\mu h^3}{12}[\sigma (\rho_{\alpha \alpha})^2+ \rho_{\alpha \beta } \rho_{\alpha \beta }]+\frac{5\mu h}{12}\varphi_{\alpha }\varphi_{\alpha }\Bigr\} \dd{a} -\int_{\Omega } \Bigl( fu+\frac{\sigma h^2}{10}f\rho_{\alpha \alpha} \Bigr) \dd{a}
\end{equation} 
among all deflections $u(x_1,x_2)$ and rotation angles $\varphi_\alpha(x_1,x_2)$ satisfying the kinematic boundary conditions. Here $\dd a=\dd x_1 \dd x_2$, $\sigma = \frac{\lambda}{\lambda +2\mu}=\frac{\nu }{1-\nu}$, with $\lambda$, $\mu$ and $\nu$ being the Lame constants and Poisson's ratio, respectively. Note that $\varphi_\alpha$ are the rotation angles due to pure shear deformation, $\rho_{\alpha \beta}$ are the measures of bending and $f$ is the external transverse load. They are given by
\begin{equation}\label{bendingmeasure}
\begin{split}
&\rho_{\alpha \beta} = u_{,\alpha \beta}-\varphi _{(\alpha ,\beta )}, \\
&f = \tau_3 |_{x_3=h/2} +\tau_3 |_{x_3=-h/2},
\end{split}
\end{equation}
with $\tau_3$ being the normal traction. In the following, we use Greek indices, running from 1 to 2, to refer to the plane coordinates $x_1$ and $x_2$.  The comma before an index denotes differentiation with respect to the corresponding coordinate, the parentheses surrounding a pair of indices denote symmetrization operation and summation over repeated indices is understood.

The FE-implementation of this variational problem is not straightforward due to the different magnitudes of the bending and shear stiffnesses. Since the shear stiffness \gb{$(5\mu h/12)$} is two orders of magnitude larger than the bending stiffness \gb{$(\mu h^3/12)$} with respect to $h$, its multiplication by very small rotation angles leads to a numerical instability called the shear-locking effect when $h$ becomes small. To solve this problem, we first try to get rid of the second derivatives in this variational problem. We introduce the new unknown functions
\begin{equation}\label{totalangle}
\psi_\alpha =-u_{,\alpha }+\varphi_\alpha
\end{equation} 
which have the meaning of the total angles of rotation of the transverse fibers due to both bending and shear deformation so that
\begin{equation}\label{kinematics}
\begin{split}
&\rho_{\alpha \beta} = u_{,\alpha \beta}-\varphi _{(\alpha ,\beta )}=-\psi_{(\alpha,\beta )}, \\
&\varphi _\alpha= u_{,\alpha }+\psi_\alpha .
\end{split}
\end{equation} 
In terms of these new unknown functions, the functional becomes
\begin{multline}\label{energyfunctional1}
J[u,\psi _\alpha ]=\int_{\Omega }\Bigl\{ \frac{\mu h^3}{12}[\sigma (\psi_{\alpha ,\alpha})^2+ \psi_{(\alpha ,\beta)} \psi_{(\alpha ,\beta)}]+\frac{5\mu h}{12}(u_{,\alpha }+\psi_\alpha )(u_{,\alpha }+\psi_\alpha )\Bigr\} \dd{a} 
\\
-\int_{\Omega } \Bigl( fu-\frac{\sigma h^2}{10}f \psi_{\alpha ,\alpha} \Bigr)\dd{a}.
\end{multline} 
To avoid the shear-locking effect, we try to make the problem independent of $h$ so that the orders of magnitude of the bending and shear stiffnesses become the same. For this purpose, we introduce the rescaled coordinates and rotation angles as follows
\begin{equation}\label{rescaling}
\bar{x}_\alpha =\frac{x_\alpha}{h},\quad \bar{\psi}_{\bar{\alpha}} =h\psi_\alpha .
\end{equation}
Note that $\bar{x}_\alpha$ are dimensionless, while $\bar{\psi}_{\bar{\alpha}} $ have the dimension of length and can be interpreted as the longitudinal displacements of the positive face surface of the plate \cite{berdichevsky1979variational,le2023asymptotically} (see Fig.~1). Since 
\begin{equation}\label{after}
\begin{split}
&\rho_{\alpha \beta}=-\psi_{(\alpha ,\beta )}=-\frac{1}{h^2}\bar{\psi}_{(\bar{\alpha },\bar{\beta})},
\\
&u_{,\alpha }+\psi_\alpha =\frac{1}{h}(u_{,\bar{\alpha} } +\bar{\psi}_{\bar{\alpha }}),
\\
&\dd{a}=\dd{x_1}\dd{x_2}=h^2\dd{\bar{x}_1}\dd{\bar{x}_2}=h^2\dd{\bar{a}},
\end{split}
\end{equation}
we reduce the energy functional to
\begin{multline}\label{energyfunctional2}
J[u,\bar{\psi }_{\bar{\alpha}} ]=h\int_{\bar{\Omega }}\Bigl\{ \frac{\mu }{12}[\sigma (\bar{\psi}_{\bar{\alpha },\bar{\alpha}})^2+ \bar{\psi}_{(\bar{\alpha },\bar{\beta})} \bar{\psi}_{(\bar{\alpha },\bar{\beta})}]+\frac{5\mu}{12}(u_{,\bar{\alpha }}+\bar{\psi}_{\bar{\alpha }})(u_{,\bar{\alpha }}+\bar{\psi}_{\bar{\alpha }})\Bigr\} \dd{\bar{a}} 
\\
-\int_{\bar{\Omega }} h^2 \Bigl( fu-\frac{\sigma }{10}f \bar{\psi}_{\bar{\alpha },\bar{\alpha}} \Bigr) \dd{\bar{a}}.
\end{multline} 
Here $\bar{\Omega}=\{ (\bar{x}_1,\bar{x}_2) \, | \, (x_1,x_2)\in \Omega\}$ denotes the rescaled 2-D domain, and components of vectors and their derivatives having Greek indices with bar are related to the rescaled coordinates $\bar{x}_\alpha $. We can further simplify this functional by dividing it by $\mu h$. Then the minimization problem reduces to
\begin{multline}\label{eq:1}
\bar{J}[u,\bar{\psi }_{\bar{\alpha }} ]=\int_{\bar{\Omega }}\Bigl\{ \frac{1}{12}[\sigma (\bar{\psi}_{\bar{\alpha },\bar{\alpha}})^2+ \bar{\psi}_{(\bar{\alpha },\bar{\beta})} \bar{\psi}_{(\bar{\alpha },\bar{\beta})}]+\frac{5}{12}(u_{,\bar{\alpha }}+\bar{\psi}_{\bar{\alpha }})(u_{,\bar{\alpha }}+\bar{\psi}_{\bar{\alpha }})\Bigr\} \dd{\bar{a}} 
\\
-\int_{\bar{\Omega }} \Bigl( \bar{f}u-\frac{\sigma }{10}\bar{f} \bar{\psi}_{\bar{\alpha },\bar{\alpha}} \Bigr) \dd{\bar{a}} \rightarrow \min_{u,\bar{\psi}_{\bar{\alpha }}},
\end{multline}
where
\begin{equation}\label{fbar}
\bar{f}=\frac{hf}{\mu }.
\end{equation} 
Thus, $\bar{f}$ is equal to the plate thickness $h$ times the characteristic strain $\varepsilon=f/\mu$, and since both the bending and shear stiffness in the rescaled functional \eqref{eq:1} have the order of unity, the minimizer must be of the same order as $\bar{f}$ multiplied by a function depending on the characteristic size of $\bar{\Omega}$. Returning to the original functions, we see that the rotation angles $\psi_\alpha$ have the order of the characteristic strain ($\varphi_\alpha$ are even much smaller \cite{le1999vibrations,berdichevsky2009variational}), while the bending measures $\rho_{\alpha \beta}$ have the order of the characteristic strain divided by $h$, as discussed in the introduction. In this way, the rescaled problem \eqref{eq:1} provides an elegant and effective way to avoid the shear-locking effect. Note that the same rescaled formulation can be applied to the FSDT for anisotropic plates, since both the kinematic formulas \eqref{totalangle}-\eqref{kinematics} and the orders of the bending and shear stiffnesses with respect to $h$ remain unchanged compared to the isotropic plates.

To solve problem \eqref{eq:1} we need to pose the boundary conditions. If the part of the plate's edge, $\bar{\partial}_k$, is clamped, the admissible functions must satisfy the following ``hard'' kinematical conditions
\begin{equation}\label{kinbc}
u=0, \quad \bar{\psi }_{\bar{\alpha }} =0 \quad \text{at $\bar{\partial }_k$}.
\end{equation}
If the remaining part of the plate's edge, $\bar{\partial}_s$, is simply supported, then only the kinematical condition $u=0$ should be fulfilled, while $\bar{\psi }_{\bar{\alpha }}$ can be varied arbitrarily.\footnote{This is the so-called ``soft'' boundary conditions for simply supported edge. About other boundary conditions for FSDT see \cite{haeggblad1990}.} Finally, if the remaining part $\bar{\partial}_s$ of the plate's edge is free, then no constraints are imposed on $u$ and $\bar{\psi }_{\bar{\alpha }}$ there. After finding the solution, we want also to check the asymptotic accuracy of the theory. For this purpose, we need to compute the true average displacement of the plate by adding the correction term to the previously determined deflection \cite{berdichevsky1979variational}
\begin{equation}\label{trueu}
\check{u}=u+\frac{h^2\sigma }{60}\rho_{\alpha \alpha}=u-\frac{\sigma }{60}\bar{\psi}_{\bar{\alpha },\bar{\alpha}}.
\end{equation}
The determination of $\varphi_\alpha$ also requires the first derivatives of $u$. Therefore we would like to have the solution of \eqref{eq:1} such that $u$ and $\bar{\psi}_{\bar{\alpha}}$ are both $C^1$-functions. To check the asymptotic accuracy we must compare functions $\check{u}$ and $\check{\psi}_\alpha$ found by the 2-D FSDT with
\begin{equation}\label{eq:ex2_displacement_rotation}
\langle w_3 (x_\alpha ,x_3)\rangle \equiv \frac{1}{h}\int_{-h/2}^{h/2} w_3 (x_\alpha ,x_3)\dd{x_3} \quad \text{and} \quad \langle w_\alpha  (x_\alpha ,x_3) x_3 \rangle /(h^2/12) ,
\end{equation}
where $w_i(x_\alpha,x_3)$ are the displacements computed by the 3-D exact theory of elasticity. If they agree with the accuracy up to $h^2/l^2$, then the asymptotic accuracy of our FE-implementation is guaranteed. 

\section{Finite element implementation}

\subsection{Weak and strong formulations}
Since we will only deal with rescaled coordinates and rotation angles in this Section and in the analytical parts of the next Section, we briefly omit all bars in the functional \eqref{eq:1}. Calculating the first variation of \eqref{eq:1} and equating it to zero, we obtain the following necessary condition for the minimum
\begin{multline}\label{eq:2}
\var{J}=\int_{\Omega }\Bigl\{ \frac{1}{6}[\sigma \psi_{\alpha ,\alpha} \var{\psi}_{\beta,\beta} + \psi_{(\alpha ,\beta)} \var{\psi}_{(\alpha ,\beta)}]+\frac{5}{6}(u_{,\alpha }+\psi_\alpha )(\var{u}_{,\alpha }+\var{\psi}_\alpha )\Bigr\} \dd{a} 
\\
-\int_{\Omega } \Bigl( f\var{u}-\frac{\sigma }{10}f \var{\psi}_{\alpha ,\alpha} \Bigr) \dd{a}=0.
\end{multline}
Introducing the notations
\begin{equation}
\label{eq:3}
\begin{split}
&\var{W}^b=\int_{\Omega }\frac{1}{6}[\sigma \psi_{\alpha ,\alpha} \var{\psi}_{\beta,\beta} + \psi_{(\alpha ,\beta)} \var{\psi}_{(\alpha ,\beta)}] \dd{a} ,
\\
&\var{W}^s=\int_{\Omega }\frac{5}{6}(u_{,\alpha }+\psi_\alpha )(\var{u}_{,\alpha }+\var{\psi}_\alpha ) \dd{a} ,
\\
&\var{W}^e= \int_{\Omega } (f\var{u}-\frac{\sigma }{10}f \var{\psi}_{\alpha ,\alpha})\dd{a},
\end{split}
\end{equation}
for the variations of the bending energy, shear energy, and the work of external forces, respectively, we recast Eq.~\eqref{eq:2} in the form
\begin{equation}
\label{eq:4}
\var{W}^b+\var{W}^s=\var{W}^e.
\end{equation}
Denote the space of admissible functions as $\mathcal{K}=\{ (v,\chi_\alpha )\, | \, (v,\chi_\alpha )|_{\partial_k}=0\}$. The weak formulation of the problem is stated as follows: Given $f$, find $(u,\psi_\alpha )\in \mathcal{K}$ such that Eq.~\eqref{eq:4} is satisfied for all $(\var{u},\var{\psi}_\alpha )\in \mathcal{K}$. For the integrals in \eqref{eq:3} to be meaningful, $v$ and $\chi_\alpha$ must belong at least to the Sobolev's space of square integrable functions with the square integrable first derivatives, $H^1(\Omega)$. However, since the desired asymptotic accuracy of FSDT may require higher smoothness of $u$ and $\psi_\alpha$, the continuity assumption in $\mathcal{K}$ still remains unspecified.

Although not relevant for the FE-implementation, the strong formulation of the problem is also stated for the sake of completeness. Under the assumptions that $u$ and $\psi_\alpha$ are doubly differentiable, we apply the partial integration to Eq.~\eqref{eq:2} and put it in the form
\begin{multline}\label{var}
\var{J}=\int_{\Omega } \Bigl\{ \Bigl[ -\frac{\sigma }{6}\psi_{\beta ,\beta \alpha} -\frac{1}{6}\psi_{(\alpha ,\beta )\beta } \Bigr] \var{\psi}_{\alpha} -\frac{5}{6}(u_{,\alpha }+\psi_\alpha )_{,\alpha }\var{u} +\frac{5}{6}(u_{,\alpha }+\psi_\alpha )\var{\psi}_\alpha \Bigr\} \dd{a}
\\
 -\int_{\Omega } \Bigl( f\var{u}+\frac{\sigma }{10}f_{,\alpha} \var{\psi}_{\alpha } \Bigr) \dd{a}+\int_{\partial_s}\Bigl[ \Bigl( \frac{\sigma }{6}\psi_{\beta ,\beta }n_\alpha +\frac{1}{6}\psi_{(\alpha ,\beta )}n_\beta \Bigr) \var{\psi}_{\alpha} 
\\
+\frac{5}{6}(u_{,\alpha }+\psi_\alpha )n_{\alpha }\var{u}\Bigr] \dd{s}+\int_{\partial_s}\frac{\sigma }{10}fn_{\alpha} \var{\psi}_{\alpha }\dd{s}=0,
\end{multline}
where $n_\alpha$ is the unit normal outward to the boundary and $\dd s$ is the length element. It is assumed here that the remaining part $\partial_s$ of the plate edge is free. Because of the arbitrariness of the variations $\var{\psi}_\alpha $ and $\var{u}$ in $\Omega$ and on $\partial_s$, from \eqref{var} follow the second order partial differential equations
\begin{equation}
\label{equil}
\begin{split}
&-\frac{\sigma +1}{6}\psi_{\beta ,\beta \alpha} -\frac{1}{6}\psi_{\alpha ,\beta \beta }  +\frac{5}{6}(u_{,\alpha }+\psi_\alpha )=\frac{\sigma }{10}f_{,\alpha},
\\
&-\frac{5}{6}(u_{,\alpha }+\psi_\alpha )_{,\alpha }=f,
\end{split}
\end{equation}
the kinematic boundary conditions \eqref{kinbc} on $\partial_k$, and the natural boundary conditions
\begin{equation}
\label{bc}
\begin{split}
&\frac{\sigma }{6}\psi_{\beta ,\beta }n_\alpha +\frac{1}{6}\psi_{(\alpha ,\beta )}n_\beta +\frac{5}{6}(u_{,\alpha }+\psi_\alpha )=-\frac{\sigma }{10}fn_{\alpha},
\\
&\frac{5}{6}(u_{,\alpha }+\psi_\alpha )n_{\alpha }=0
\end{split}
\end{equation}
on $\partial_s$. Eqs.~\eqref{equil}, \eqref{kinbc}, and \eqref{bc} constitute the strong formulation of the problem.

\subsection{Discretization}

To discretize Eq.~\eqref{eq:4}, interpolation spaces for displacement and rotation angles are required. In the following, a general interpolation scheme is used without regard to specific requirements on the order and continuity of the function space as mentioned above. In this sense, $u$ and $\psi_\alpha $ ($\alpha =1,2$) are considered as primary variables. Their interpolation is
\begin{equation}
\label{eq:5}
u=(\vb{N}^u)^T \vb{u}, \quad \psi_\alpha =(\vb{N}^\psi)^T \vb*{\psi}_\alpha,
\end{equation}
where $\vb{N}^u$ and $\vb{N}^\psi$ are column vectors representing the shape functions for  $u$ and $\psi_\alpha $, respectively, while $\vb{u}$ and $\vb*{\psi}_\alpha$ are deflection and rotation angle vectors, respectively. The vector notation is in column ordering unless otherwise stated.

Based on Eq.~\eqref{eq:5} we find the variation of the primary variables
\begin{equation}
\label{eq:6}
\var{u}=\var{\vb{u}}^T\vb{N}^u , \quad \var{\psi}_\alpha =\var{\vb*{\psi}}^T_\alpha \vb{N}^\psi ,
\end{equation}
and their derivatives
\begin{equation}\label{derivatives}
\begin{split}
&\var{u}_{,\alpha}=\var{\vb{u}}^T\vb{N}^u_{,\alpha},
\\
&\var{\psi}_{\alpha,\alpha}= \var{\vb*{\psi}}^T_\alpha \vb{N}^\psi_{,\alpha }=\var{\vb*{\psi}}^T_1 \vb{N}^\psi_{,1} + \var{\vb*{\psi}}^T_2 \vb{N}^\psi_{,2},
\\
&\var{\psi}_{(\alpha,\beta)}= \frac{1}{2}(\var{\vb*{\psi}}^T_\alpha \vb{N}^\psi_{,\beta }+\var{\vb*{\psi}}^T_\beta \vb{N}^\psi_{,\alpha}).
\end{split}
\end{equation}
Substituting Eqs.~\eqref{eq:6} into Eqs.~\eqref{eq:3}, we obtain
\begin{equation}\label{energyvariation}
\begin{split}
\var{W}^b&=\frac{1}{6}\int_{\Omega }[\sigma (\var{\vb*{\psi}}^T_1 \vb{N}^\psi_{,1} + \var{\vb*{\psi}}^T_2 \vb{N}^\psi_{,2}) (\psi_{1,1}+\psi_{2,2})  + \var{\vb*{\psi}}^T_1 \vb{N}^\psi_{,1} \psi_{1,1} 
\\
&+\var{\vb*{\psi}}^T_2 \vb{N}^\psi_{,2}\psi_{2,2}+(\var{\vb*{\psi}}^T_1 \vb{N}^\psi_{,2}+\var{\vb*{\psi}}^T_2 \vb{N}^\psi_{,1})
\psi_{(1,2)}] \dd{a} ,
\\
\var{W}^s&=\frac{5}{6}\int_{\Omega }[(\var{\vb{u}}^T\vb{N}^u_{,1}+\var{\vb*{\psi}}^T_1 \vb{N}^\psi ) (u_{,1}+\psi_1) 
\\
&+(\var{\vb{u}}^T\vb{N}^u_{,2}+\var{\vb*{\psi}}^T_2 \vb{N}^\psi ) (u_{,2}+\psi_2)]\dd{a} ,
\\
\var{W}^e&= \int_{\Omega } [\var{\vb{u}}^T\vb{N}^u f -\frac{\sigma }{10} (\var{\vb*{\psi}}^T_1 \vb{N}^\psi_{,1} + \var{\vb*{\psi}}^T_2 \vb{N}^\psi_{,2})f] \dd{a}.
\end{split}
\end{equation}
Consequently, the residual forces with respect to $u$ and $\psi_\alpha$ read
\begin{equation}
\label{eq:9}
\begin{split}
\vb{R}^u&=\int_{\Omega }\vb{N}^u f\dd{a}  -\frac{5}{6}\int_{\Omega }[\vb{N}^u_{,1}(u_{,1}+\psi_1)  + \vb{N}^u_{,2}(u_{,2}+\psi_2) ]\dd{a},
\\
\vb{R}^{\psi_1}&=-\int_{\Omega }\frac{\sigma }{10} \vb{N}^\psi_{,1}f \dd{a} -\frac{5}{6}\int_{\Omega } \vb{N}^\psi (u_{,1}+\psi_1) \dd{a}
\\
&-\frac{1}{6}\int_{\Omega }[ \sigma \vb{N}^\psi_{,1} (\psi_{1,1}+\psi_{2,2})+\vb{N}^\psi_{,1} \psi_{1,1}+\vb{N}^\psi_{,2} \psi_{(1,2)}] \dd{a},
\\
\vb{R}^{\psi_2}&=-\int_{\Omega }\frac{\sigma }{10} \vb{N}^\psi_{,2}f \dd{a} -\frac{5}{6}\int_{\Omega } \vb{N}^\psi (u_{,2}+\psi_2) \dd{a}
\\
&-\frac{1}{6}\int_{\Omega }[ \sigma \vb{N}^\psi_{,2} (\psi_{1,1}+\psi_{2,2})+\vb{N}^\psi_{,2} \psi_{2,2}+\vb{N}^\psi_{,1} \psi_{(1,2)}] \dd{a}.
\end{split}
\end{equation}

Taking the directional derivatives of Eqs.~\eqref{eq:9}, we derive the following blocks of the stiffness matrix
\begin{equation}\label{block1}
\begin{split}
\vb{K}^{uu}&= \frac{5}{6}\int_{\Omega } \left[ \vb{N}^u_{,1}(\vb{N}^u_{,1})^T + \vb{N}^u_{,2}(\vb{N}^u_{,2})^T \right] \dd{a},
\\
\vb{K}^{u \psi_1}&=\frac{5}{6}\int_{\Omega } \vb{N}^\psi_{,1}(\vb{N}^\psi)^T \dd{a},
\\
\vb{K}^{u \psi_2}&=\frac{5}{6}\int_{\Omega } \vb{N}^\psi_{,2}(\vb{N}^\psi)^T \dd{a},
\end{split}
\end{equation}
\begin{equation}\label{block2}
\begin{split}
\vb{K}^{\psi_1 u}&= \frac{5}{6}\int_{\Omega }\vb{N}^\psi (\vb{N}^\psi_{,1})^T\dd{a},
\\
\vb{K}^{\psi_1 \psi_1}&=\int_{\Omega } \Bigl[ \frac{5}{6}\vb{N}^\psi (\vb{N}^\psi)^T +\frac{\sigma+1}{6} \vb{N}^\psi_{,1} (\vb{N}^\psi_{,1})^T 
+\frac{1}{12} \vb{N}^\psi_{,2} (\vb{N}^\psi_{,2})^T\Bigr] \dd{a},
\\
\vb{K}^{\psi_1 \psi_2}&=\int_{\Omega } \Bigl[ \frac{\sigma}{6} \vb{N}^\psi_{,1} (\vb{N}^\psi_{,2})^T +\frac{1}{12} \vb{N}^\psi_{,2} (\vb{N}^\psi_{,1})^T\Bigr] \dd{a},
\end{split}
\end{equation}
\begin{equation}\label{block3}
\begin{split}
\vb{K}^{\psi_2 u}&= \frac{5}{6}\int_{\Omega }\vb{N}^\psi (\vb{N}^\psi_{,2})^T\dd{a},
\\
\vb{K}^{\psi_2 \psi_1}&=\int_{\Omega } \Bigl[ \frac{\sigma}{6} \vb{N}^\psi_{,2} (\vb{N}^\psi_{,1})^T 
+\frac{1}{12} \vb{N}^\psi_{,1} (\vb{N}^\psi_{,2})^T\Bigr] \dd{a},
\\
\vb{K}^{\psi_2 \psi_2}&=\int_{\Omega } \Bigl[ \frac{5}{6}\vb{N}^\psi (\vb{N}^\psi)^T +\frac{\sigma+1}{6}  \vb{N}^\psi_{,2} (\vb{N}^\psi_{,2})^T 
+\frac{1}{12} \vb{N}^\psi_{,1} (\vb{N}^\psi_{,1})^T\Bigr] \dd{a}.
\end{split}
\end{equation}
The final form of the discretized linear system reads
\begin{equation}
\label{final}
\begin{pmatrix}
  \vb{K}^{uu}   & \vb{K}^{u \psi_1}  &  \vb{K}^{u \psi_2} \\
 \vb{K}^{\psi_1 u}  &  \vb{K}^{\psi_1 \psi_1}  &  \vb{K}^{\psi_1 \psi_2}  \\
 \vb{K}^{\psi_2 u}  &  \vb{K}^{\psi_2 \psi_1}  &  \vb{K}^{\psi_2 \psi_2}
\end{pmatrix} \begin{pmatrix}
 \vb{u}    \\
\vb*{\psi}_1  \\
\vb*{\psi}_2   
\end{pmatrix}
=
\begin{pmatrix}
\vb{F}^u    \\
\vb{F}^{\psi_1}  \\
\vb{F}^{\psi_2}  
\end{pmatrix},
\end{equation}
with the elements of the stiffness matrix being the blocks given in Eqs.~\eqref{block1}-\eqref{block3}, and $\vb{F}^u$, $\vb{F}^{\psi_1}$, $\vb{F}^{\psi_2}$ the first terms in \eqref{eq:9} representing the external forces. 

\subsection{Isogeometric analysis}

The weak form \eqref{var} contains the first order derivatives of deflection and rotation angles and their variations and is therefore meaningful only if $(u,\psi_\alpha)$ and $(\var{u},\var{\psi}_\alpha)$ belong to the Sobolev space $H^1(\Omega)$. However, as mentioned in Section 2, for the desired asymptotic accuracy, the discretization space must be at least $C^1$ to ensure continuity and smoothness of the solution. In this paper, we use isogeometric analysis with the non-uniform rational
B-splines (NURBS) shape function, which guarantees $C^1$-continuity \cite{Hughes_Cottrell_Bazilevs:05}. This approach also allows to facilitate the refinement of the mesh and the elevation of the discretization order if necessary.

In the context of geometric modeling, a plate can be represented by NURBS surface patches, in which each patch can be described as a tensor product of two univariate B-splines such as
\begin{equation}
\mathbf{S} \left( \xi_1, \xi_2 \right) = \sum_{i=1}^m \sum_{j=1}^n N_i^p \left( \xi_1 \right) N_j^q \left( \xi_2 \right) \mathbf{P}_{ij}.
\label{eq:nurbs_surface}
\end{equation}
In Eq.~\eqref{eq:nurbs_surface}, the univariate B-spline basis function of order $p$ and $q$, denoted as $N_i^p$ and $N_j^q$, is computed via recursive Cox-de-Boor formula
\begin{equation}\label{boor}
N_i^p(\xi) = \dfrac{\xi-\xi_i}{\xi_{i+p}-\xi_{i}} N^{p-1}_i(\xi) + \dfrac{\xi_{i+p+1}-\xi_i}{\xi_{i+p+1}-\xi_{i+1}} N^{p-1}_{i+1}(\xi) , \; N_i^0(\xi)=\begin{cases}
1 & \xi_i \leq \xi \leq \xi_{i+1} \\
0 & \text{otherwise}
\end{cases}.
\end{equation}
We use $\{ \mathbf{P}_{ij} \}_{0 \leq i \leq m, 0 \leq j \leq n}$ to denote the control point grid. The definition of the control point in homogeneous coordinates will facilitate the construction of the NURBS patch. A surface can comprise a single or multiple patches that are connected at the patch interface. In the latter case, the concept is called multipatch analysis. In the typical scenario, the patch information matches at the interface, including the parametric (knot) information and the location of the control points, the multipatch is naturally connected and the $C^0$ continuity is guaranteed on the patch interfaces, i.e. strong coupling. However, patch-wise parametric matching is not the necessary condition. There are various methods to maintain the patch continuity in a NURBS multipatch, such as the penalty method \cite{Marussig.Hughes:17}, Nitsche method \cite{Ruess.etal:14} or bending strip method \cite{Kiendl.etal:10}. These methods impose weak coupling condition and are more sophisticated to implement.

The design of the B-spline surface allows different interpolation orders in each parametric direction. This feature is particularly useful for analyses of materials with anisotropic behavior where accuracy needs to be improved in a particular direction. In addition, rescaling of the mesh can be conveniently performed by rescaling the control point coordinates. This eliminates the need for remeshing.

The NURBS patches in a multipatch structure are macro elements where the shape function needs the full knot vectors for values and derivatives evaluation. This on one hand prevents the possibility of applying parallelization algorithms to the numerical code, on the other hand it makes the element matrices very dense, which reduces the performance of the sparse direct solver. B\'ezier decomposition algorithm \cite{Borden.etal:11} is proposed to retain the local characteristic of the finite element. In essence, the B\'ezier decomposition strategy constructs a local operator for each knot span (in 1D) or knot cell (in 2D and 3D) so that the shape function within the cell can be evaluated independently with the knot vectors. The support domain, i.e. control points influenced by the cell, of each cell and the number of shape functions required for each cell is reduced to $p+1$ in each parametric direction. This eliminates the global access to patch-wise data and allows a typical finite element code to work seamlessly with NURBS elements.

\begin{figure}[!htb]
\centering
\includegraphics[scale=1]{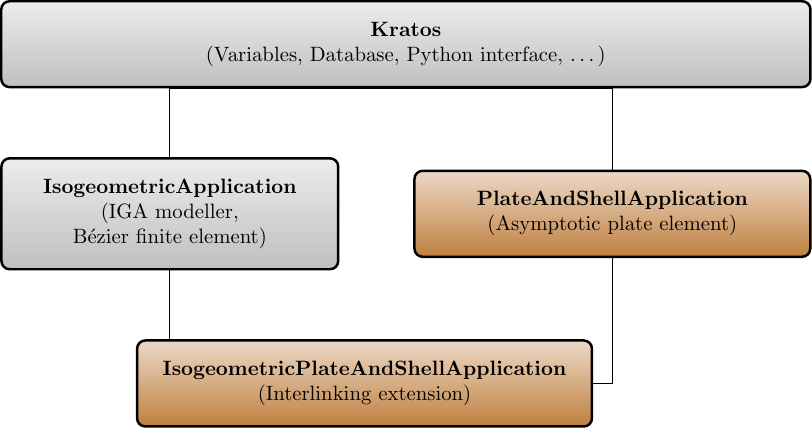}
\caption{Structure of the computational code.}
\label{fig:program_structure}
\end{figure}

\subsection{Structure of the finite element code}

The computational code for simulating the proposed FSDT plate element is implemented as an extension of the Kratos Multiphysics framework \cite{Dadvand.etal:10}. A brief overview of the structure of the code is visualized in Fig.~\ref{fig:program_structure}. The plate element is implemented within the extension \texttt{PlateAndShellApplication}. By combining with the existing extension \texttt{IsogeometricApplication}, which supports NURBS patch modeling and B\'ezier extraction operator evaluation, a new extension \texttt{IsogeometricPlateAndShellApplication} can be created. This module enables isogeometric analysis with FSDT plate element. In addition, mesh refinement and order elevation are supported via \texttt{IsogeometricApplication}.

\section{Numerical examples}

\subsection{Circular plate under uniform transverse load}

As a first test problem, we apply the developed finite element code to the numerical simulation of the axisymmetric bending of a circular plate clamped at the edge and subjected to a constant transverse load $f$. We chose this 2-D problem because it admits an analytical solution, so we can directly compare the numerical solution with it to test the former's divergence and investigate the optimal choice of interpolation order. As the bending is axisymmetric, we must have in the polar coordinates $(r,\theta)$ the vanishing rotation angle $\psi_\theta=0$, while $u$ and $\psi_r$ are functions of $r$ only. Consequently,
\begin{equation}\label{rho}
\rho_{rr}=-\psi_{r,r},\quad \rho_{\theta \theta}=-\frac{1}{r}\psi_r,
\end{equation}
while all other components of $\rho_{\alpha \beta}$ vanish. It is easy to show that the problem \eqref{eq:1} reduces to minimizing the following functional
\begin{multline}\label{energyfunc3}
I[u,\psi _r]=\int_{0}^R\Bigl\{ \frac{1}{12}[\sigma (\psi_{r,r}+\frac{1}{r}\psi_r)^2+ (\psi_{r,r})^2+\frac{1}{r^2}(\psi_{r})^2] +\frac{5}{12}(u_{,r}+\psi_{r})^2\Bigr\} r \dd{r} 
\\
-\int_{0}^R \Bigl[ f u-\frac{\sigma }{10}f \Bigl( \psi_{r,r}+\frac{1}{r}\psi_r \Bigr) \Bigr] r\dd{r} \rightarrow \min_{u(R)=\psi_{r}(R)=0} ,
\end{multline}
where $R$ is the dimensionless radius of the plate and all bars are still dropped for short. The standard calculus of variations shows that the minimizer must satisfy the equation
\begin{equation}
\label{eq:15}
\begin{split}
& -\frac{\sigma}{6}\Bigl[ (\psi_{r,r}+\frac{1}{r}\psi_r)r\Bigr]_{,r}+\frac{\sigma}{6}(\psi_{r,r}+\frac{1}{r}\psi_r)-\frac{1}{6}(\psi_{r,r}r)_{,r}
\\
&+\frac{1}{6r}\psi_r+\frac{5}{6}(u_{,r}+\psi_{r})r=0,
\\
&-\frac{5}{6}\left[ (u_{,r}+\psi_{r})r \right]_{,r}=fr,
\end{split}
\end{equation} 
subjected to the clamped boundary conditions $u(R)=\psi_{r}(R)=0$. Integrating the second equation of \eqref{eq:15} and using the fact that $\varphi=u_{,r}+\psi_{r}$ cannot be singular at $r=0$, we find that
\begin{equation}
\label{eq:16}
\varphi(r)=u_{,r}+\psi_{r}=-\frac{3}{5}fr.
\end{equation}
Substituting \eqref{eq:16} into the first of \eqref{eq:15}, we reduce it to
\begin{equation}\label{reduced}
\psi_{r,rr}+\frac{1}{r}\psi_{r,r}-\frac{1}{r^2}\psi_r=-\frac{3}{\sigma+1}fr.
\end{equation}
Integrating this equation and using the boundary condition $\psi_r(R)=0$ as well as the non-singularity of $\psi$ at $r=0$, we get
\begin{equation}
\label{eq:18}
\psi_r(r)=\frac{3}{8(\sigma+1)}f(R^2r-r^3)=\frac{3(1-\nu)}{8}f(R^2r-r^3).
\end{equation}
Now, the deflection $u$ can be found from 
\begin{equation}
u_{,r}=-\frac{3}{5}fr-\psi_r=-\frac{3}{5}fr-\frac{3}{8(\sigma+1)}f(R^2r-r^3).
\end{equation}
Integrating this equation and using the boundary condition $u(R)=0$, we obtain finally (cf. \cite{wang1996deflection,reddy2022theories})
\begin{equation}\label{deflection}
u(r)=\frac{3}{10}f(R^2-r^2)+\frac{3}{32(\sigma+1)}f(R^2-r^2)^2=\frac{3}{10}f(R^2-r^2)+\frac{3(1-\nu)}{32}f(R^2-r^2)^2.
\end{equation}
Note that the deflection according to Kirchhoff's plate theory reads \cite{landau1986theory}
\begin{equation}
u_K(r)=\frac{3(1-\nu)}{32}f(R^2-r^2)^2.
\end{equation}
while the angle of rotation, $-u_{K,r}$, coincides with \eqref{eq:18}.

If we have the simply supported edge instead of the clamped edge, then the deflection becomes
\begin{equation}
\label{supported}
u(r)=\frac{3}{10(1+\nu)}f(R^2-r^2)+\frac{3(1-\nu)}{32}f(R^2-r^2)\Bigl( \frac{5+\nu}{1+\nu}R^2-r^2\Bigr) .
\end{equation} 
Since the solution procedure is similar to the previous case, we omit it here for brevity. Note that the second term on the right-hand side of \eqref{supported} represents the deflection according to Kirchhoff's plate theory \cite{landau1986theory}.

\begin{figure}[!htb]
\centering
\includegraphics[scale=0.14,trim=20cm 2cm 19cm 1cm,clip]{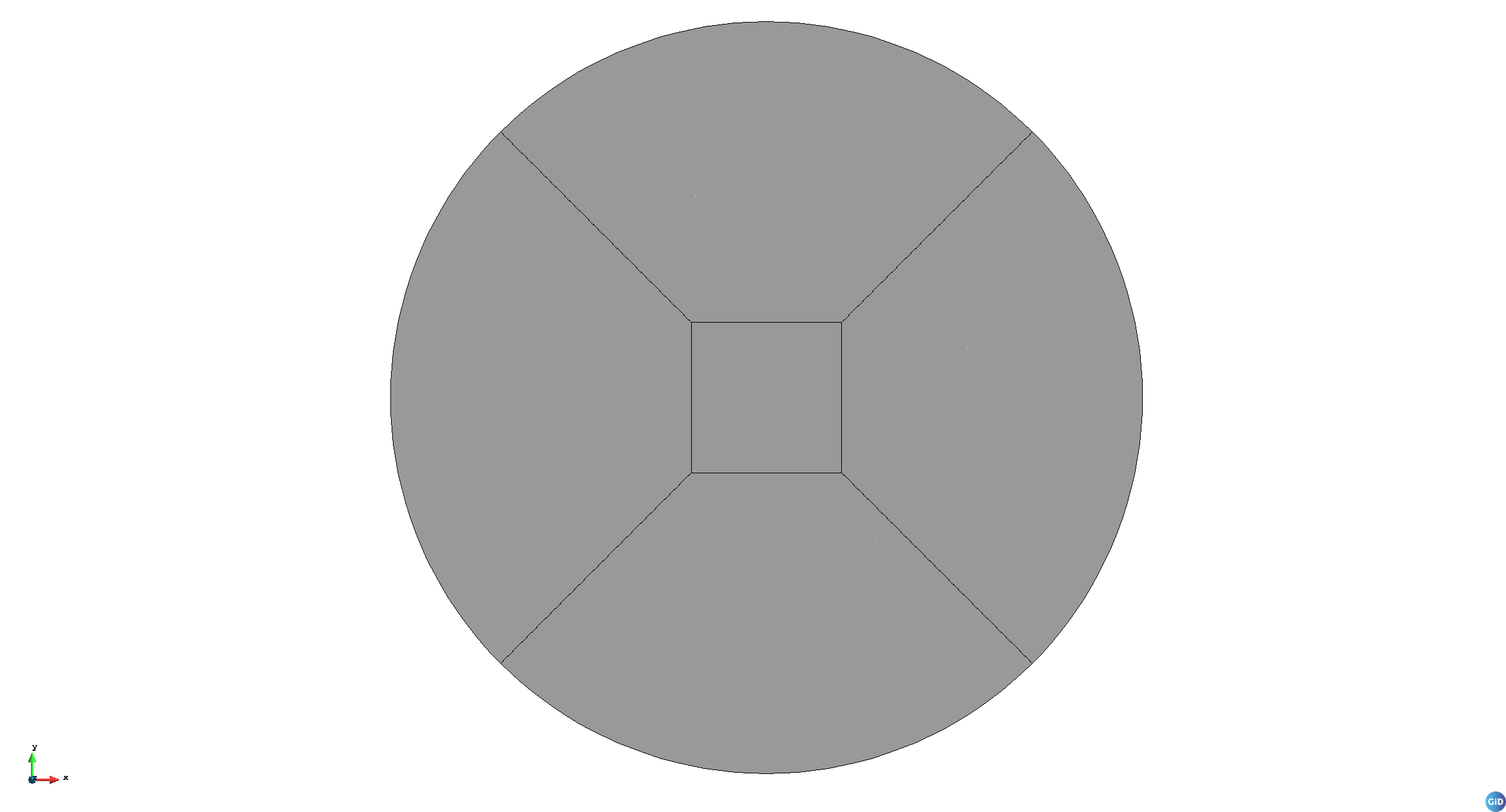}
\includegraphics[scale=0.14,trim=20cm 2cm 19cm 1cm,clip]{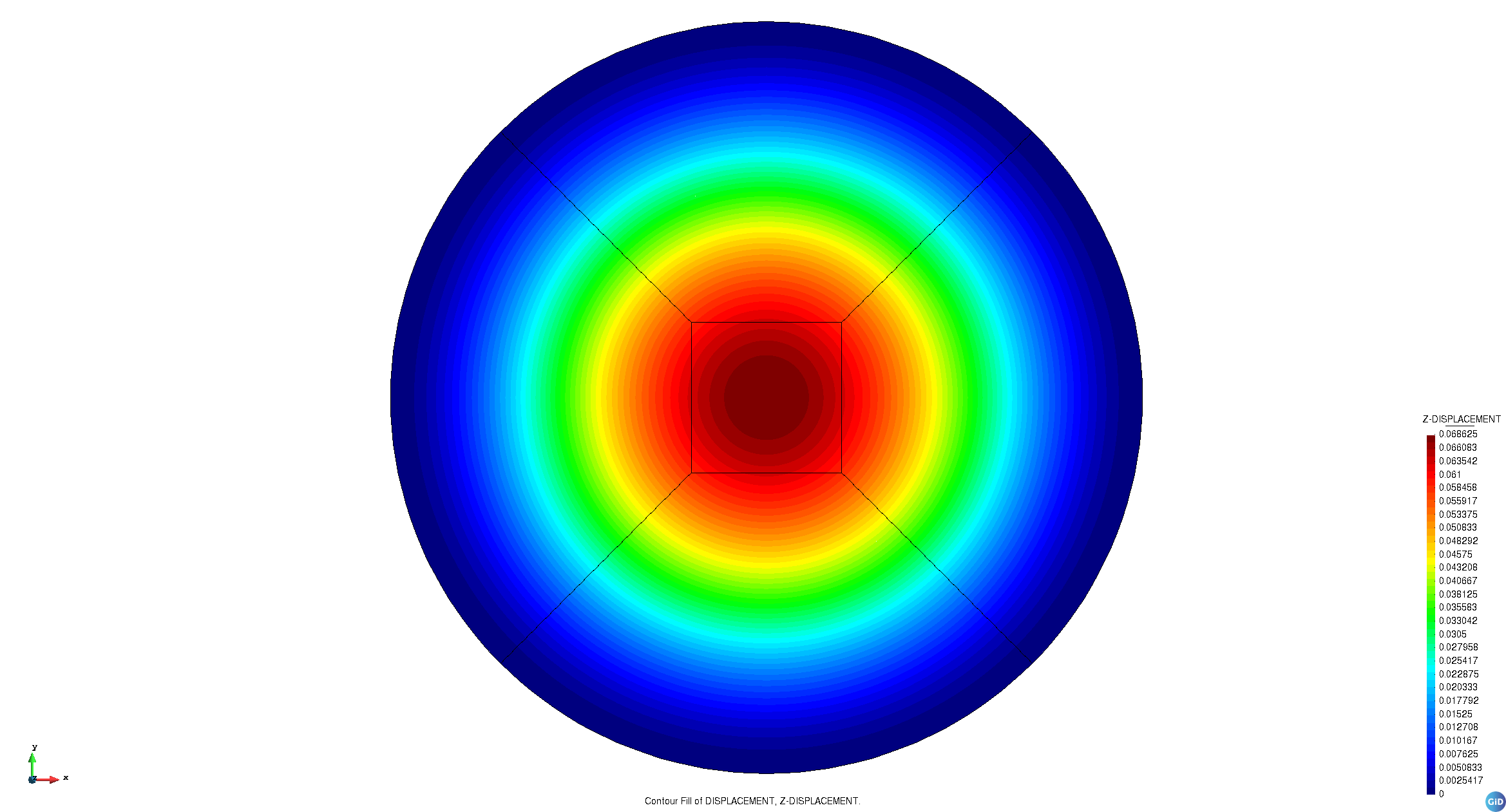}
\includegraphics[scale=0.3,trim=74cm 2cm 1cm 22.5cm,clip]{disp_z_scaled.png}
\caption{Numerical example 4.1(a) (Circular plate with clamped edge): Multipatch geometry (left) and contour plot of scaled deflection $\left( \dfrac{u}{h \varepsilon \bar{R}^4} \right)$ at $\bar{R}=10$ (right).}
\label{fig:ex1_geometry}
\end{figure}

For the numerical simulations of the real plate let us return to the original notation with bars. In the rescaled formulation, the problem for the circular plate contains only three parameters: $\bar{R}=R/h$ (geometry), $\nu$ (material), and $\bar{f}=hf/\mu$ (load). The material parameter is chosen to be $\nu=0.3$. Since the problem is linear, the solution depends linearly on $\bar{f}$, so by setting $\bar{f}=1$ we calculate the normalized solution. To get the real deflection, we just need to multiply our normalized $u$ by $h\varepsilon$ (where $\varepsilon=f/\mu$). Similarly, to calculate the real rotation angles, we need to multiply our normalized $\bar{\psi}_{\bar{\alpha}}$ (or $\bar{\varphi}_{\bar{\alpha}}$) by $\varepsilon$. To investigate the efficiency of the proposed plate element, we perform the analysis with three dimensionless radii $\bar{R}=10$, $\bar{R}=100$ and $\bar{R}=1000$. The circular plate is constructed with 5 NURBS fields, as shown in Fig.~\ref{fig:ex1_geometry} (left). 

\begin{figure}[!htb]
\centering
\includegraphics[scale=1]{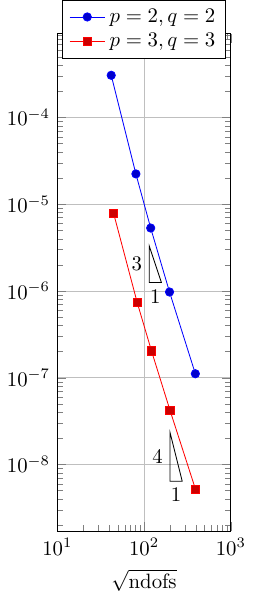}
\includegraphics[scale=1]{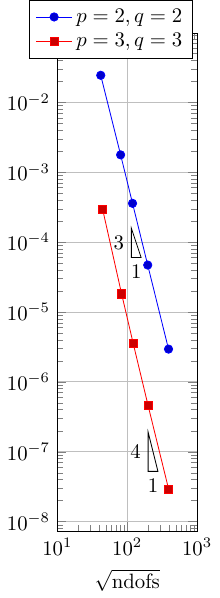}
\includegraphics[scale=1]{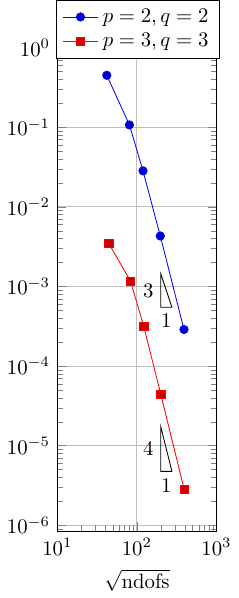}
\caption{Numerical example 4.1(a) (Circular plate with clamped edge): Convergence rate of displacement ($L_2$ error) with $\bar{R}=10$ (left), $\bar{R}=100$ (middle) and $\bar{R}=1000$ (right).}
\label{fig:ex1_convergence_rate}
\end{figure}

The deflection agrees very well with the analytical solution \eqref{deflection} (with bars being restored over the rescaled quantities), as shown in Table~\ref{tab:ex1_w_results}. The numerical values of the normalized deflection (for $\bar{f}=1$) are large. However, considering that the real deflection is obtained after multiplication by $h\varepsilon$, where $\varepsilon$ is much smaller than $10^{-3}$, these values become reasonable. The $L_2$ error shown in the right column of this Table is computed according to
\begin{equation}
\| e \|_{L_2} = \dfrac{\int_{\Omega} \left( u_{computed} - u_{ana} \right)^2 \dd{a}}{\int_{\Omega} u^2_{ana} \dd{a}}.
\end{equation}
It can be seen that the cubic order discretization gives very good results even at coarse mesh. Fig.~\ref{fig:ex1_convergence_rate} visualizes the convergence rate of analysis with second order ($p=2,q=2$) and cubic order ($p=3,q=3$) discretization. We observe that the error increases with increasing $\bar{R}$, yet the convergence rate is close to optimal, showing that the proposed formulation does not exhibit shear-locking behavior. Although rescaling inflates the domain, the number of elements and degrees of freedom required to obtain a small error are not significantly affected by this rescaling. A slight reduction of convergence rate is observed for $\bar{R}=10$ of the cubic order analysis. This can be attributed to the loss of accuracy when the $L_2$ error is very small ($\sim 10^{-9}$) and the machine precision starts to affect the computation.

\begin{table}[!htb]
\centering
\begin{tabular}{cccccccc}
\hline
\multicolumn{5}{c}{$\bar{R}=10$} \\
\cline{2-5}
& ndofs & $u$ & $u_e$ \% & $L_2$ error \\
\hline
\multirow{5}{*}{IGA-p2-q2}
& 1752  &  6.8615046723e2 & 1.4503863949e-2 &  3.0660060196e-4 \\
& 6492  &  6.8624375111e2 & 9.1058452571e-4 & 2.2378738174e-5 \\
& 14232  & 6.8624876504e2 & 1.7995811726e-4 & 5.3214286086e-6 \\
& 38712  & 6.8624983992e2 & 2.3326748357e-5 & 9.7389087961e-7 \\
& 152412 & 6.8624999000e2 & 1.4579201519e-6 & 1.1120739899e-7 \\
\hline
\multirow{5}{*}{IGA-p3-q3}
& 1986  &  6.8624841430e2 & 2.3106803895e-4 &  7.8773270723e-6 \\
& 6936  &  6.8624990465e2 & 1.3894075265e-05 &  7.4005944177e-7 \\
& 14886 &  6.8624998130e2 & 2.7245754934e-06 &  2.0315875429e-7 \\
& 39786  & 6.8624999759e2 & 3.5171608593e-07 &  4.1964948985e-8 \\
& 154536 & 6.8624999985e2 & 2.1741422667e-08 &  5.1396533350e-9 \\
\hline
\multicolumn{5}{c}{$\bar{R}=100$} \\
\cline{2-5}
& ndofs & $u$ & $u_e$ & $L_2$ error \\
\hline
\multirow{5}{*}{IGA-p2-q2}
& 1752 &   6.4784590440e6  & 1.3257323276 & 2.4336817628e-2 \\
& 6492  &  6.5593384992e6  & 9.3846634028e-02 & 1.7699590666e-3 \\
& 14232  & 6.5642574481e6  & 1.8925472721e-02 & 3.5821555521e-4 \\
& 38712  & 6.5653375422e6  & 2.4744165221e-03 & 4.6899163890e-5 \\
& 152412 & 6.5654898199e6  & 1.5505456947e-04 & 2.9411995117e-6 \\
\hline
\multirow{5}{*}{IGA-p3-q3}
& 1986  &  6.5644455940e6 & 1.6059797387e-2 &  2.9372761135e-4 \\
& 6936  &  6.5654300922e6 & 1.0647752621e-03 &  1.8101151749e-5 \\
& 14886 &  6.5654861376e6 & 2.1113940168e-04 &  3.5340802949e-6 \\
& 39786 &  6.5654982029e6 & 2.7372178611e-05 &  4.5584241877e-7 \\
& 154536 & 6.5654998877e6 & 1.7097563813e-06 &  2.8764590547e-8 \\
\hline
\multicolumn{5}{c}{$\bar{R}=1000$} \\
\cline{2-5}
& ndofs & $u$ & $u_e$ & $L_2$ error \\
\hline
\multirow{5}{*}{IGA-p2-q2}
& 1752  &  4.6439563968e10 & 2.9235273640e1 &  4.4941755443e-1 \\
& 6492  &  6.1501718808e10 & 6.2835235678e+00 &  1.0699473685e-1 \\
& 14232 &  6.4595537398e10 & 1.5691548868e+00 &  2.8321455202e-2 \\
& 38712  & 6.5473621571e10 & 2.3112797835e-01 &  4.3069485094e-3 \\
& 152412 & 6.5615277512e10 & 1.5272292485e-02 &  2.8832105167e-4 \\
\hline
\multirow{5}{*}{IGA-p3-q3}
& 1986  &  6.5493593646e10 & 2.0069447938e-1 &  3.4437565289e-3 \\
& 6936  &  6.5582249763e10 & 6.5600061147e-02 &  1.1614533849e-3 \\
& 14886  & 6.5613643904e10 & 1.7761589296e-02 &  3.1440539208e-4 \\
& 39786  & 6.5623595741e10 & 2.5969542594e-03 &  4.4297675346e-5 \\
& 154536 & 6.5625189222e10 & 1.6880430048e-04 &  2.8107143827e-6 \\
\end{tabular}
\caption{Numerical example 4.1(a) (Circular plate with clamped edge): Deflection computed at the center of the plate and error. $u_e$ is the error of deflection compared to the analytical solution.}
\label{tab:ex1_w_results}
\end{table}

\begin{figure}[!htb]
\centering
\includegraphics[scale=0.55]{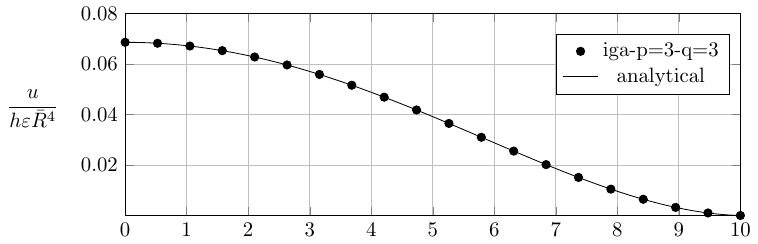}
\includegraphics[scale=0.55]{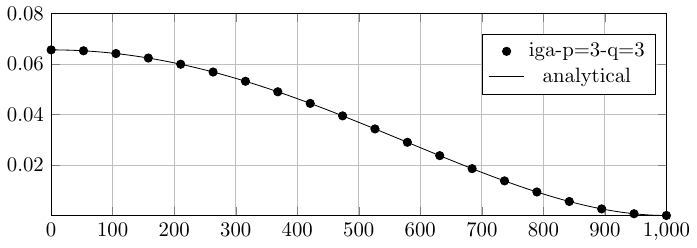} \\
\includegraphics[scale=0.55]{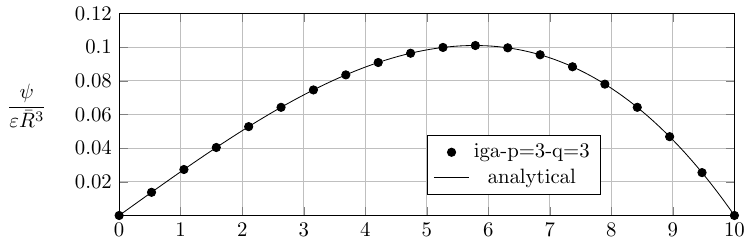}
\includegraphics[scale=0.55]{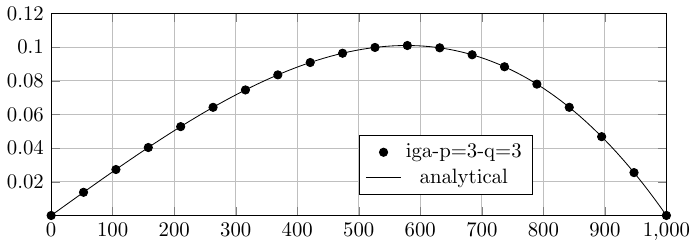} \\
\includegraphics[scale=0.55]{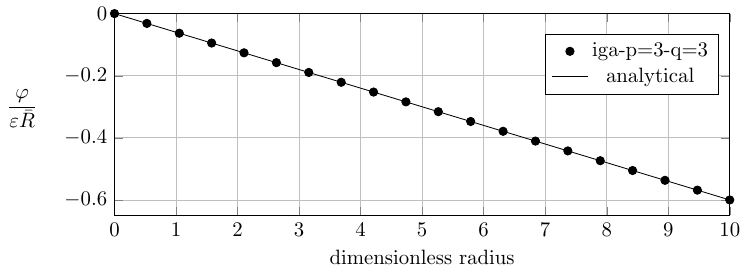}
\includegraphics[scale=0.55]{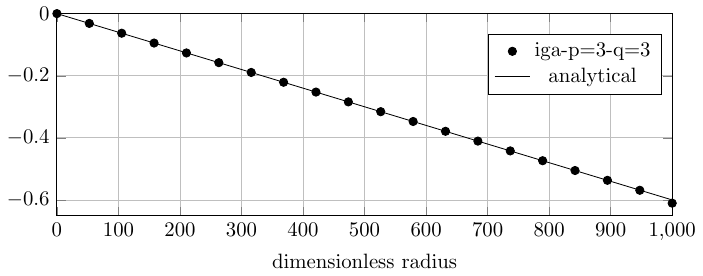}
\caption{Numerical example 4.1(a) (Circular plate with clamped edge): Normalized deflection and rotation angle along the radial direction for $\bar{R}=10$ (left) and $\bar{R}=1000$ (right).}
\label{fig:ex1_rotation}
\end{figure}

Due to the high accuracy of the solution, the distribution of the deflection is almost radially symmetric as expected (see Fig.~\ref{fig:ex1_geometry} (right) for the solution of the cubic discretization with 39786 d.o.fs for $\bar{R}=10$). The normalized deflection $u/(h \varepsilon \bar{R}^4)$ together with the normalized rotation angles, i.e. $\psi_{r}/(\varepsilon \bar{R}^3)$ and $\varphi_{r}/(\varepsilon \bar{R})$, for this particular discretization along the radial direction are plotted in Fig.~\ref{fig:ex1_rotation}. It clearly shows perfect agreement with the analytical formulas \eqref{eq:16}-\eqref{deflection} (where the bars are recovered over the scaled quantities) for both $\bar{R}=10$ and $\bar{R}=1000$. It is worth noting that to evaluate $\varphi$, the derivatives of $u$ at the integration point are calculated and then transferred to the control points using an $L_2$-global projection algorithm \cite{Jiao.Heath:04}.

\begin{figure}[!htb]
\centering
\includegraphics[scale=1]{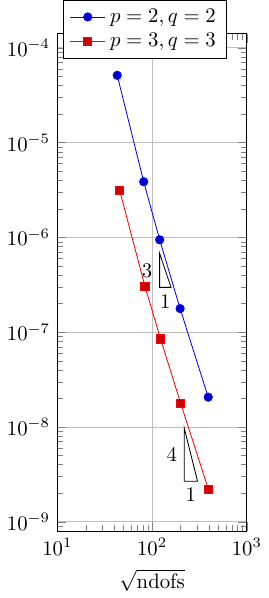}
\includegraphics[scale=1]{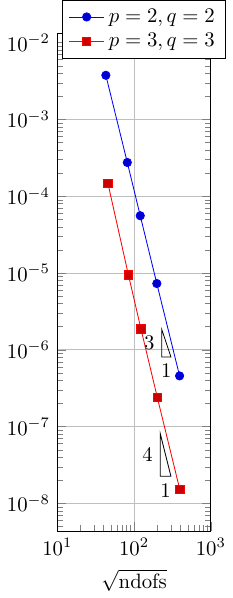}
\includegraphics[scale=1]{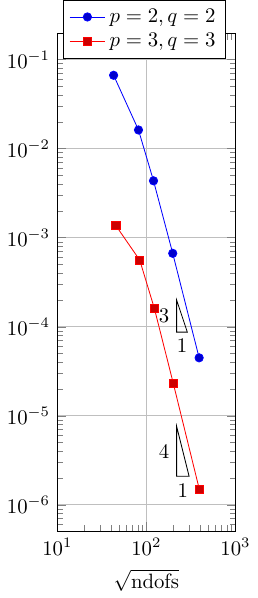}
\caption{Numerical example 4.1(b) (Circular plate with simply supported edge): Convergence rate of displacement ($L_2$ error) with $\bar{R}=10$ (left), $\bar{R}=100$ (middle) and $\bar{R}=1000$ (right).}
\label{fig:ex1_ss_convergence_rate}
\end{figure}

Fig.~\ref{fig:ex1_ss_convergence_rate} shows the convergence rate to the analytical solution described by \eqref{supported} in the analysis of the circular plate with simply supported edge. We observe the same optimal convergence behavior as in Fig.~\ref{fig:ex1_convergence_rate}. 

\subsection{Rectangular plate under uniform transverse load}
\begin{figure}[!htb]
\centering
\includegraphics[scale=0.8]{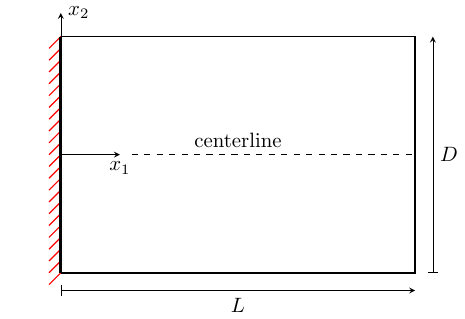}
\raisebox{10mm}[0pt][0pt]{%
\includegraphics[scale=0.91]{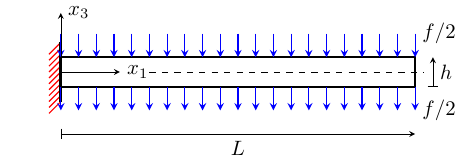}
}
\caption{Numerical example 4.2(a) (Rectangular plate with one clamped edge and three free edges): Geometry and boundary condition for (left) the analysis using FSDT plate element and (right) the 2D analysis using solid element.}
\label{fig:ex2_beam}
\end{figure}

For the second test problem, we apply the developed finite element code to the numerical simulation of the bending of a rectangular plate subjected to a constant load $f$. We assume that $\Omega$ occupies the domain $(0,L)\times (-D/2,D/2)$ of the $(x_1,x_2)$-plane. In the first case study we let the left edge at $x_1=0$ be clamped while the three remaining edges of the plate be free (see Fig.~\ref{fig:ex2_beam} (left)).

When the depth of the rectangular plate, $D$, becomes large, the solution must exhibit the plane strain state due to the almost translational invariance in the $x_2$-direction. For FSDT in the rescaled formulation \eqref{eq:1}, this means that the 2-D problem reduces to the following 1-D variational problem: Minimize the energy functional of the beam-like model
\begin{multline}\label{eq:21}
I[u,\psi ]=\int_{0}^L\Bigl[ \frac{1}{2}E_b (\psi_{,x})^2+ \frac{1}{2}E_s(u_{,x}+\psi)^2\Bigr] \dd{x} 
\\
-\int_{0}^L \Bigl( f u-\frac{\sigma }{10}f \psi_{,x}\Bigr) \dd{x} \rightarrow \min_{u(0)=\psi(0)=0} ,
\end{multline}
where the bars in this theoretical part \gb{are} again dropped for short, $x\equiv x_1$ and
\begin{equation}
E_b=\frac{\sigma+1}{6}=\frac{1}{6(1-\nu)},\quad E_s=\frac{5}{6}.
\end{equation}
The vanishing first variation leads to the differential equations
\begin{equation}\label{eq:22}
\begin{split}
& -E_b \psi_{,xx}+E_s(u_{,x}+\psi)=0,
\\
&-E_s (u_{,x}+\psi)_{,x}=0,
\end{split}
\end{equation} 
which are subjected to the boundary conditions
\begin{equation}
\begin{split}
& u(0)=\psi(0)=0,
\\
&E_b \psi_{,x}(L)=-\frac{\sigma}{10}f,\quad E_s(u_{,x}+\psi)|_{x=L}=0.
\end{split}
\end{equation} 
From the first equation of \eqref{eq:22} follows
\begin{equation}
u_{,x}+\psi=\varphi=\frac{E_b}{E_s}\psi_{,xx}.
\end{equation}
Substituting this into the second equation yield
\begin{equation}
\label{eq:23}
E_b\psi_{,xxx}=-f.
\end{equation}
Integrating Eq.~\eqref{eq:23} and using the boundary conditions, we obtain
\begin{equation}\label{eq:24}
\psi(x)=\frac{f}{E_b}\Bigl[ -\frac{1}{6}x^3+\frac{1}{2}Lx^2-(\frac{1}{2}L^2+\frac{\sigma}{10})x \Bigr] .
\end{equation}
The deflection can be found from the equation
\begin{equation}
u_{,x}=\varphi-\psi=-\psi+\frac{E_b}{E_s}\psi_{,xx},
\end{equation}
which, together with the boundary condition $u(0)=0$, yields
\begin{equation}
u(x)=\frac{f}{E_b}\Bigl[ \frac{1}{24}x^4-\frac{1}{6}Lx^3+(\frac{1}{4}L^2+\frac{\sigma}{20})x^2 \Bigr] +\frac{f}{E_s}(-\frac{x^2}{2}+Lx).
\end{equation}
However, this is still not the true average transverse displacement of the plate. The latter should be computed in accordance with Eq.~\eqref{trueu} giving
\begin{multline}\label{eq:25}
\check{u}(x)=\frac{f}{E_b}\Bigl[ \frac{1}{24}x^4-\frac{1}{6}Lx^3+(\frac{1}{4}L^2+\frac{\sigma}{20})x^2 \Bigr] +\frac{f}{E_s}(-\frac{x^2}{2}+Lx)
\\
-\frac{\sigma}{60}\frac{f}{E_b}\Bigl(-\frac{1}{2}x^2+Lx-\frac{1}{2}L^2-\frac{\sigma}{10}\Bigr) .
\end{multline}
Note that the average transverse displacement according to the classical Kirchhoff's plate theory in this beam-like model is given by
\begin{equation}
u_{K}(x)=\frac{f}{E_b}\Bigl( \frac{1}{24}x^4-\frac{1}{6}Lx^3+\frac{1}{4}L^2x^2 \Bigr) ,
\end{equation} 
while the angle of rotation reads
\begin{equation}
\psi_{K}(x)=-u_{K,x}=-\frac{f}{E_b}\Bigl( \frac{1}{6}x^3-\frac{1}{2}Lx^2+\frac{1}{2}L^2x \Bigr) .
\end{equation}

The presence of the analytical solution given by Eqs.~\eqref{eq:24}-\eqref{eq:25} for this particular case allows a comparison with the solution of the 2-D FSDT and with the solution of the exact elasticity theory. On one hand, the comparison with the solution of the 2-D FSDT can show that our FE-code is free of shear-locking and allows us to check its convergence and efficiency as $D$ becomes large. On the other hand, the comparison with the solution of the exact elasticity theory allows us to verify the asymptotic accuracy of our FE-code. To solve this problem within the exact elasticity theory, we use the geometry and boundary condition shown in Fig.~\ref{fig:ex2_beam} (right). The normal traction is applied to the top and bottom of the plate, as shown in this Figure, to avoid its elongation and at the same time maintain the bending combined with a pure shear \cite{berdichevsky1979variational}. Under the plane strain condition the problem becomes 2-D, so we can employ 2-D solid elements. Having the numerical solution within the plane strain elasticity theory, the average transverse displacement and rotation angles over the thickness are evaluated using Eq.~\eqref{eq:ex2_displacement_rotation}. Within 2-D FSDT the deflection and rotation angles are evaluated along the centerline $x_2=0$. Then the true average transverse displacement $\check{u}$ is computed in accordance with \eqref{trueu}.

\begin{figure}[!htb]
\centering
\includegraphics[scale=0.57]{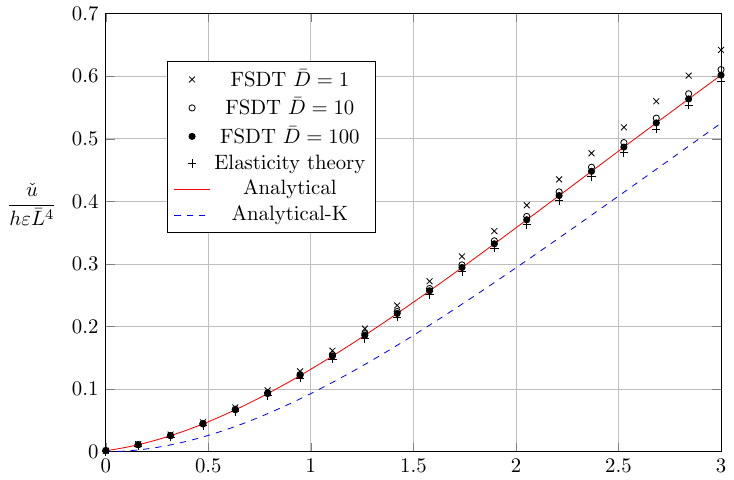}
\includegraphics[scale=0.57]{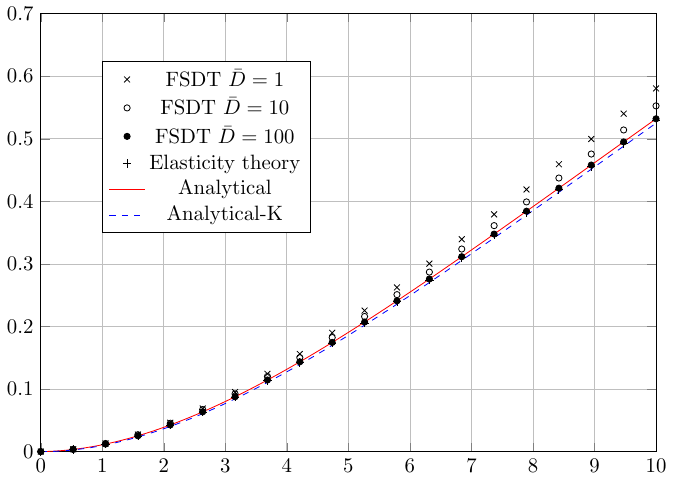}
\includegraphics[scale=0.55]{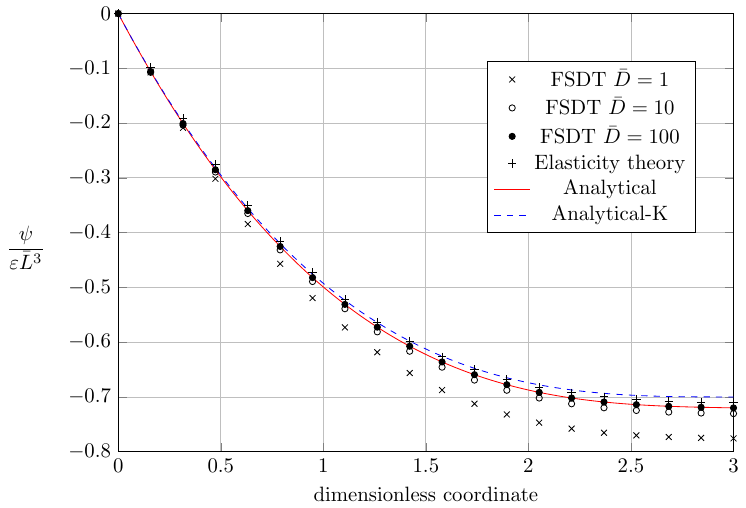}
\includegraphics[scale=0.55]{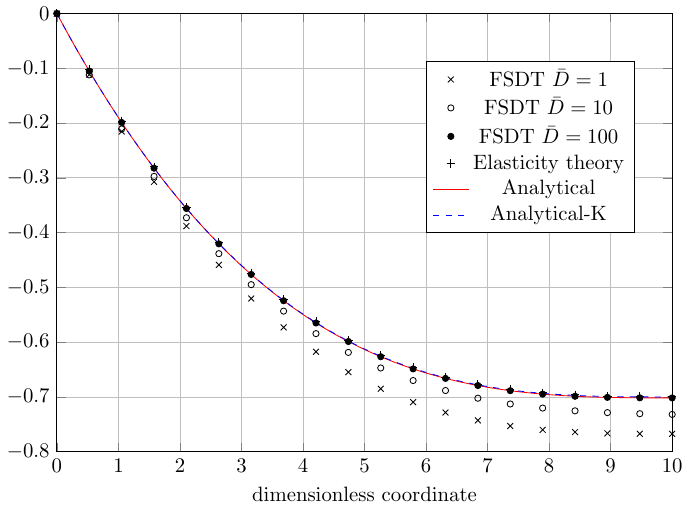}
\caption{Numerical example 4.2(a) (Rectangular plate with one clamped edge and three free edges): Average displacement and rotation angle for  $\bar{L}=3$ (left) and $\bar{L}=10$ (right).}
\label{fig:ex2_L3_u_psi}
\end{figure}

Both the plate analysis and the elastic solid analysis employ the cubic order NURBS discretization with sufficiently fine mesh. Using again the bars for the recalled quantities, we show the results of numerical simulation in Fig.~\ref{fig:ex2_L3_u_psi}, in which two cases $\bar{L}=3$ (left) and $\bar{L}=10$ (right) are analyzed (with $\bar{L}=L/h$).

In Fig.~\ref{fig:ex2_L3_u_psi}, two trends are clearly visible. First, the numerical solution using the 2-D FSDT plate element converges to the analytical solution of FSDT given by Eqs.~\eqref{eq:24}-\eqref{eq:25} (with the bars being recovered over the scaled quantities) as $\bar{D}=D/h$ increases, as shown in all plots for three different cases $\bar{D}=1$, $\bar{D}=10$, and $\bar{D}=100$. This is due to the fact that the boundary effect becomes negligible for large $\bar{D}$ and agrees well with the assumption of almost translational invariance discussed at the beginning of this sub-Section. This also confirms that our rescaled formulation of the FSDT is inherently shear-locking-free and that the FE-code is efficient. Second, the solution of the 2-D FSDT approximates the solution of the exact elasticity theory much better than that of the classical Kirchhoff plate theory. For large $\bar{L}$ (thin plate), the difference between the solutions of FSDT and exact elasticity theory is negligibly small. Even for moderate $\bar{L}$ (moderately thick plate), the error is only noticeable near the end point when the deflection reaches its maximum. This also confirms the asymptotic accuracy of our FE-implementation of the FSDT. Since the difference between the analytical solution of the FSDT (red solid line) and the solution resulting from the classical Kirchhoff (beam-like) theory (blue dashed line) becomes significant for moderate $\bar{L}$, the latter is not applicable for moderately thick plates. A detail worth mentioning is that the true average displacement calculated according to \eqref{trueu} does not exactly satisfy the kinematic boundary condition: The last correction term in \eqref{trueu}, although small, does not vanish at $\bar{x}_1=0$ (cf. Eq.~\eqref{eq:25}). The reason for this is well known: The FSDT correction might not work in a thin boundary layer near the edge of the plate \cite{berdichevsky1979variational}. Far from the edge, this correction term of order $h^2/l^2$ compared to unity is essential for ensuring the asymptotic accuracy of the FSDT up to that order. Another detail worth mentioning is that the elastic 2-D analysis does not use an approximation, but evaluates the solution fields in a brute force method with a very fine mesh. The average transverse displacements and rotation angles obtained by calculating integrals \eqref{eq:ex2_displacement_rotation} over the thickness in a post-processing represent only the integral characteristics of the FSDT. The comparison between the detailed distributions of displacements across the thickness is left for future investigations.

\begin{figure}[!htb]
\centering
\includegraphics[scale=0.9]{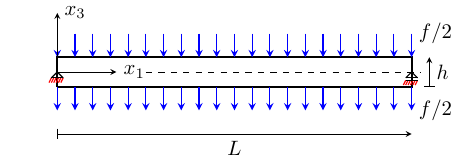}
\caption{Numerical example 4.2(b): A plate with two simply supported edges and two free edges.}
\label{fig:4.2b}
\end{figure}

In the second case study we assume that both the left and right edges of the plate be simply supported, while the two remaining edges be free (see Fig.~\ref{fig:4.2b}). For sufficiently large $D$ the 2-D problem becomes 1-D due to the same reason as in the previous case. This 1-D beam-like model admits again an analytical solution. Omitting for brevity the detailed solution procedure which is similar to the previous case, we summarize the final formulas below
\begin{multline}
\label{deflection2}
\check{u}(x)=\frac{f}{E_b}\Bigl[ \frac{1}{24}x^4-\frac{1}{12}Lx^3+\frac{\sigma}{20}x^2+(\frac{1}{24}L^3-\frac{\sigma}{20}L) x \Bigr] +\frac{f}{E_s}(-\frac{x^2}{2}+\frac{1}{2}Lx)
\\
-\frac{\sigma}{60}\frac{f}{E_b}\Bigl(-\frac{1}{2}x^2+\frac{1}{2}Lx-\frac{\sigma}{10}\Bigr) ,
\end{multline} 
\begin{equation}\label{psi2}
\psi(x)=\frac{f}{E_b}\Bigl[ -\frac{1}{6}x^3+\frac{1}{4}Lx^2-\frac{\sigma}{10}x -\frac{1}{24}L^3+\frac{\sigma}{20}L \Bigr] ,
\end{equation}
\begin{equation}\label{varphi2}
\varphi(x)=\frac{f}{E_s}\Bigl( -x +\frac{1}{2}L \Bigr) .
\end{equation}
Note that the average transverse displacement according to the classical Kirchhoff's plate theory in this beam-like model is given by
\begin{equation}
u_{K}(x)=\frac{1}{24}\frac{f}{E_b} (x^4-2Lx^3+L^3x) ,
\end{equation} 
while the angle of rotation reads
\begin{equation}
\psi_{K}(x)=-u_{K,x}=-\frac{1}{24}\frac{f}{E_b}( 4x^3-6Lx^2+L^3) .
\end{equation}

\begin{figure}[!htb]
\centering
\includegraphics[scale=0.53]{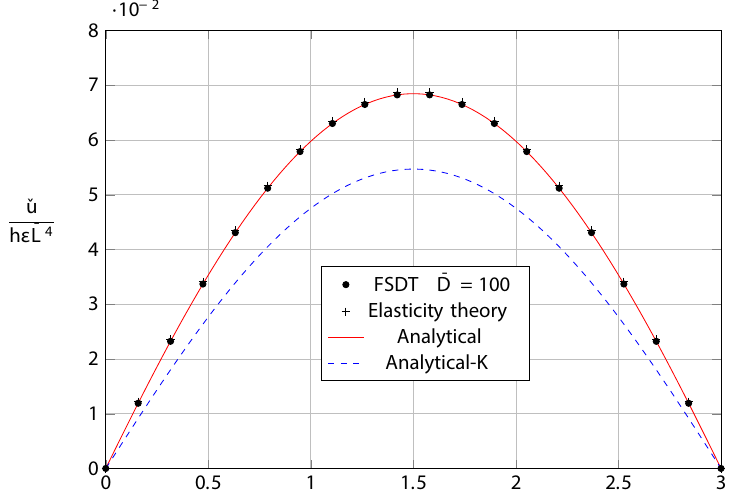}
\includegraphics[scale=0.53]{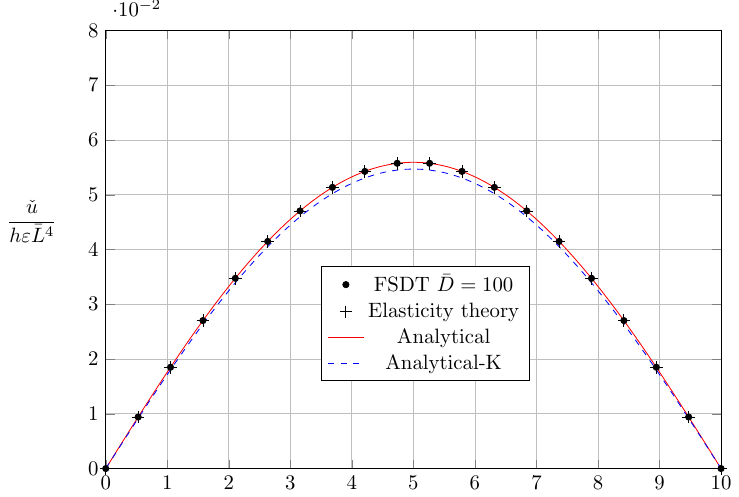}
\includegraphics[scale=0.53]{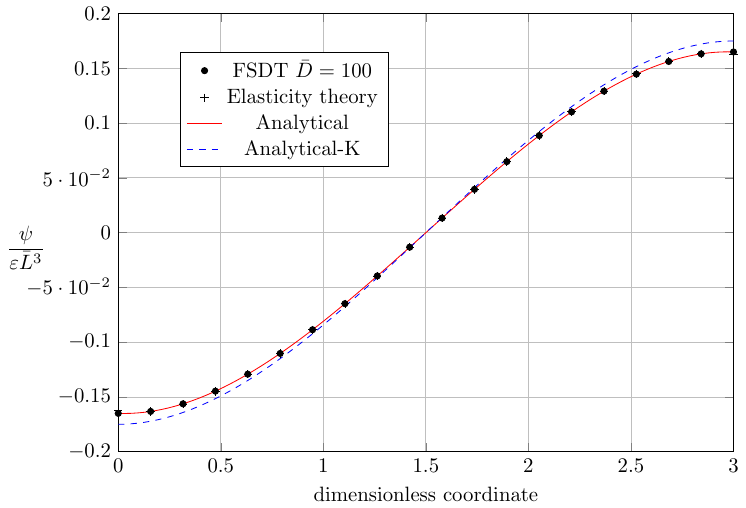}
\includegraphics[scale=0.53]{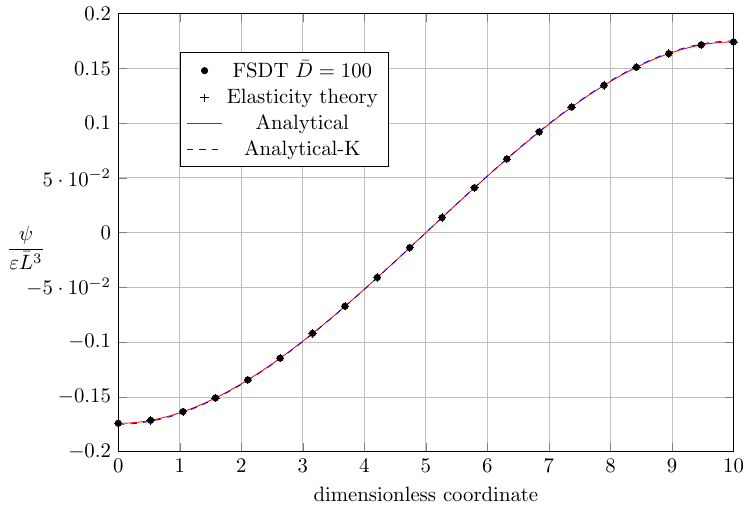}
\caption{Numerical example 4.2(b) (Rectangular plate with simply supported boundary condition): Average displacement and rotation angle for  $\bar{L}=3$ (left) and $\bar{L}=10$ (right).}
\label{fig:ex2_ss_L3_u_psi}
\end{figure}

The numerical results for this case study are shown in Fig.~\ref{fig:ex2_ss_L3_u_psi}. It can be seen that the solution of the 2-D FSDT agrees very well with the analytical solution of the 1-D beam-like model (for large $\bar{D}=100$) and with the solution found by the elasticity theory, even for moderate $\bar{L}=3$. For $\bar{L}=10$ all three theories provide almost identical results. For $\bar{L}=3$, the error of the FSDT in determining the average transverse displacement is only noticeable in the center of the plate, where the deflection reaches its maximum. Note that Kirchhoff's plate theory becomes unsuitable for moderate $\bar{L}$, as the upper left plot in Fig.~\ref{fig:ex2_ss_L3_u_psi} shows. It should also be noted that the solution determined by the 2-D analysis using solid element satisfies the following integral constraints at $x_1=0$ and $x_1=L$
\begin{equation}
\label{eq:integralbc}
\int_{-h/2}^{h/2} w_1 \dd{x_3}=0,\quad \int_{-h/2}^{h/2} w_3 \dd{x_3}=0,
\end{equation}
which reproduce the simply supported boundary conditions in elasticity theory quite well. 

\begin{figure}[!htb]
\centering
\includegraphics[scale=0.163,trim=0 6cm 5cm 0,clip]{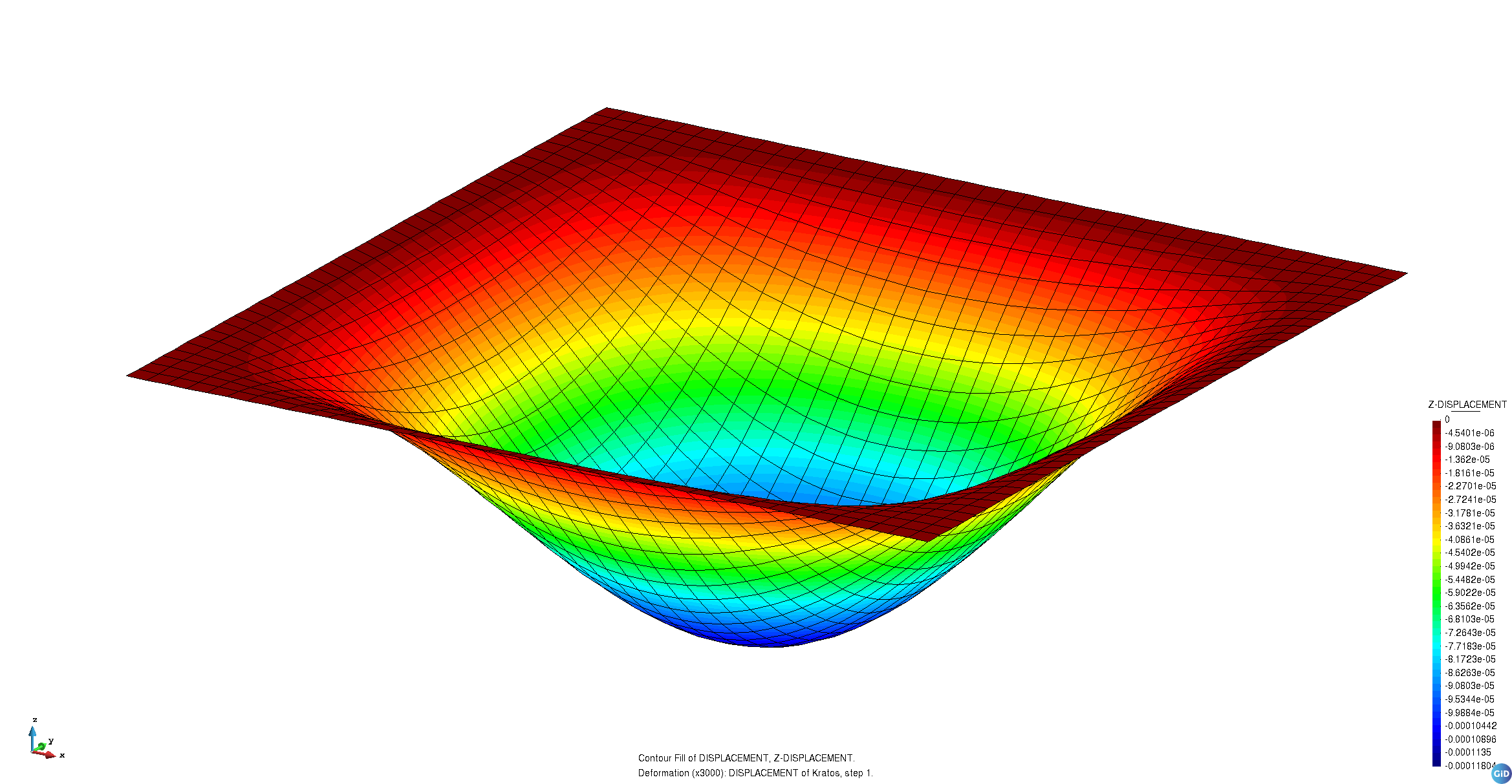}
\includegraphics[scale=0.3,trim=74.5cm 2cm 0 0,clip]{disp_z.png}
\caption{Numerical example 4.2(c) (Square concrete plate with clamped boundary condition): Deflection of the plate.}
\label{fig:ex3}
\end{figure}

Finally, for a realistic example, we perform an analysis for a square concrete plate that is clamped at four edges and deforms under its own weight. The material parameters of the concrete are $E=22.95 $GPa, $\nu=0.3$, length $L=1$m, thickness $h=5$cm and density $\rho = 2400$kg/m$^3$. The deflection is shown in Fig.~\ref{fig:ex3}. The maximum vertical deflection, measured in the center of the plate, is $\sim 0.12$mm. Note that the FE-code works for any geometry, material parameters and boundary conditions of the plate.

\section{Conclusion}
In this work, two advancements were made compared to the available approaches: (i) the inherently shear-locking free formulation of the asymptotically exact FSDT for homogeneous elastic plates in the rescaled coordinates and rotation angles was found, (ii) the asymptotic accuracy of the isogeometric FE-analysis that is computationally efficient and free of high-order interpolation schemes and/or sophisticated integration techniques has been achieved. The requirement that the deflection and rotation angles belong to the $C^1$ function space is only necessary to achieve asymptotic accuracy but is not relevant to the shear-locking effect, which is absent in this formulation. Since our main focus was on the inherently shear-locking-free formulation and the asymptotic accuracy of its FE-implementation, the simplest linear FSDT for homogeneous elastic plates was chosen intentionally. Once understood, the results obtained in this paper can be extended in many directions. Among others, we will develop FE-implementation based on the inherently shear-locking-free formulation of the asymptotically exact: (i) linear FSDT for homogeneous elastic shells, (ii) linear FSDT for FG-plates and shells, (iii) dynamic FSDT for FG-plates and shells, (iv) nonlinear FSDT for FG-plates and shells with application to buckling analysis. The asymptotically accurate and shear-locking free finite element implementation of the linear FSDT for homogeneous elastic shells will be addressed in our next publication.

\end{document}